\documentclass[a4paper, 11pt]{amsart}
\usepackage{amssymb,amsmath,txfonts,mathrsfs,titletoc,appendix, mathtools}
\usepackage[alphabetic]{amsrefs}
\usepackage{geometry}

\newtheorem{theorem}{Theorem}[section]
\newtheorem{prop}[theorem]{Proposition}
\newtheorem{lemma}[theorem]{Lemma}
\newtheorem{remark}[theorem]{Remark}

\newtheorem{question}[theorem]{Question}
\newtheorem{definition}[theorem]{Definition}
\newtheorem{cor}[theorem]{Corollary}

\numberwithin{equation}{section}

\def\pf{{\it Proof:}~}

\begin{document}

\title[An equation linking $\mathscr{W}$-entropy with reduced volume]{An equation linking $\mathscr{W}$-entropy with reduced volume}
\author{Guoyi Xu}
\address{Mathematical Sciences Center\\Jin Chun Yuan West Building \\ Tsinghua University, Beijing\\P. R. China, 100084}
\email{gyxu@math.tsinghua.edu.cn}
\date{\today}

\begin{abstract}\label{abstract}
$\mathscr{W}$-entropy and reduced volume for the Ricci flow were introduced by Perelman, which had proved their importance in the study of the Ricci flow. L. Ni studied the analogous concepts for the linear heat equation on the static manifolds, and established an equation which links the large time behavior of these two. Due to the surprising similarity between those concepts in the Ricci flow and the linear heat equation, a natural question whether such equation holds for the Ricci flow ancient solution was asked by L. Ni. In this paper, we gave an alternative proof to L. Ni's original equation based on a new method. And following the same philosophy of this method, we answer L. Ni's question positively for Type I $\kappa$-solutions of the Ricci flow. 
\\[3mm]
\end{abstract}

\maketitle

\titlecontents{section}[0em]{}{\hspace{.5em}}{}{\titlerule*[1pc]{.}\contentspage}
\titlecontents{subsection}[1.5em]{}{\hspace{.5em}}{}{\titlerule*[1pc]{.}\contentspage}
\tableofcontents

\section{Introduction}\label{section 1}
In $1982$, R. Hamilton \cite{RT} introduced the Ricci flow on compact Riemannian manifolds, which is the solution of the evolution equation deforming the metric on any $n$-dimensional Riemannian manifold $(M^n, g_{ij})$:
\begin{equation}\label{1.1.0}
{\frac{\partial}{\partial t}g_{ij}(x,t)= -2R_{ij}(x,t) 
}
\end{equation}
where $R_{ij}$ is the Ricci curvature of $M^n$ at time $t$, and it is called an ancient solution if it is defined for all $(-\infty, 0]$. In this paper, we focus on the study of ancient solution to the Ricci flow, for notation convenience, we recalled the concept of backward Ricci flow as the following:
\begin{equation}\label{1.1.1}
{\frac{\partial}{\partial t}g_{ij}(x,t)= 2R_{ij}(x,t) 
}
\end{equation}
Then an ancient solution to the Ricci flow corresponds to a backward Ricci flow solution defined on $[0, \infty)$.

To fix notation, $d(x, y, t)$ denotes the distance function with respect to the metric $g(t)$, $d \mu_{g(t)}(x)$ denotes the volume element of $g(t)$ at $x$, $B(x, r, t)$ denotes the geodesic ball of radius $r$ centered $x$ with respect to $g(t)$, and $V_{g(t)}(A)$ is the volume of the set $A\subset M^n$ with respect to the volume element $d \mu_{g(t)}$.

\begin{definition}\label{def 1.2}
{A complete, non-flat backward Ricci flow solution $(M^n, g(t))_{t\in [0, \infty)}$ is a \textbf{$\kappa$-solution} if it is \textbf{$\kappa$-noncollapsed at all scales} for some positive constant $\kappa$, i.e., for any $(x_0, t_0)\in M^n\times [0, \infty)$ and $r> 0$, if 
\begin{align}
|Rm(x, t)|\leq r^2\ , \quad \quad \quad \forall \ (x, t)\in P(x_0, t_0, r, r^{2}) \nonumber
\end{align}
where $P(x_0, t_0, r, r^{2})\vcentcolon = \{(y, s)|\ d(y, x_0, s)< r, \ t_0< s< t_0+ r^2\}$, then $V_{g(t_0)}\Big(B\big(x_0, r, t_0\big)\Big)\geq \kappa r^n$.
}
\end{definition}

In $2002$, Perelman introduced $\mathscr{W}$-entropy and reduced volume for the Ricci flow in the celebrated paper \cite{Pere}, which turn out to be of fundamental importance in the study of Ricci flow, especially for the study of the $\kappa$-solution to the Ricci flow.

Motivated by the entropy introduced by Perelman, in \cite{Ni1}, \cite{Ni2} and \cite{Ni3}, L. Ni studied the similar entropy for the linear heat equation and the reduced volume for the static metric. More concretely, let $(M^n, g)$ be a complete Riemannian manifold with $Rc\geq 0$ and maximum volume growth, namely $\theta_{\infty}\vcentcolon= \lim_{r\rightarrow \infty} \frac{V\big(B(x,r)\big)}{\omega_n r^n}> 0$, where $\omega_n$ is the volume of the unit ball of $\mathbb{R}^n$ and $V\big(B(x,r)\big)$ denotes the volume of the geodesic ball $B(x, r)$. L. Ni defined Nash entropy and $\mathscr{W}$-entropy for linear heat equation as the following:
\begin{align}
N(H, t)\vcentcolon= \int_{M^n} \Big( fH \Big) d\mu - \frac{n}{2}\ , \quad \mathscr{W}(f, t)\vcentcolon= \int_{M^n} \Big(t|\nabla f|^2+ f- n\Big)H d\mu \label{1.1.NW}
\end{align}
where $H(x, y, t)= (4\pi t)^{-\frac{n}{2}} e^{-f(x, y, t)}$ is the heat kernel on $(M^n, g)$. The monotonicity of $N(H, t)$ and $\mathscr{W}(f, t)$ was proved, and the following equation was established:
\begin{align}
\lim_{t\rightarrow \infty}\mathscr{W}(f, t)= \lim_{t\rightarrow \infty} N(H, t)= \ln \theta_{\infty} \label{1.1}
\end{align} 
From Lemma $8.10$ of \cite{RFG}, $\theta_{\infty}$ is the limit of the reduced volume for the static metric $g$,
\begin{align}
\theta_{\infty}= \lim_{t\rightarrow \infty} \int_{M^n} (4\pi t)^{-\frac{n}{2}} e^{-\frac{d^2(x, y)}{4t}} dx \label{1.1.2}
\end{align}

To prove (\ref{1.1}), L. Ni used the sharp pointwise bounds for the heat kernel proved by Li, Tam and Wang (see \cite{LTW}), which is closely related to the large time behavior of heat kernel studied by Li in \cite{Li}. 

On the other side, as observed by T. H. Colding, $\lim_{t\rightarrow \infty} t\frac{\partial}{\partial t} N(H, t)= 0$ can follow from the cone structure at infinity of
the manifold (see Section $5$ of \cite{Colding}). Also note $\mathscr{W}(f, t)= t\frac{\partial}{\partial t}N(H, t)+ N(H, t)$, this motivated our alternate proof of (\ref{1.1}) in Section \ref{Section 2} following Colding's observation. 

The key point of our proof of (\ref{1.1}) is that we only need the uniform (not sharp) Gaussian bounds of heat kernel, which were proved by Li and Yau in \cite{LY}. In our case, uniform means the coefficients in the bounds do not depend on time $t$. Using such uniform bounds of heat kernel, we can reduce the proof of (\ref{1.1}) to the tangent cone at infinity of manifold by sort of ``Dominated Convergence Theorem". From the work of Cheeger-Colding the tangent cone is a metric cone $C(X)$ (see  \cite{CC}), and the explicit formula of heat kernel on $C(X)$ was given by Y. Ding in \cite{Ding}, (\ref{1.1}) on $C(X)$ follows from these facts explicitly.  

For the Ricci flow, instead of the heat kernel to the linear heat equation, we consider the conjugate heat kernel, which is the defined as the following.
\begin{definition}\label{def chk}
{If $(M^n, g(t))_{t\in [0, \infty)}$ is a solution to the backward Ricci flow, the \textbf{conjugate heat kernel} $H(x, y, t)$ is the fundamental solution to 
\begin{equation}\label{2.2}
{\frac{\partial}{\partial t} H(x, y, t)= \Delta_x H(x, y, t)- R(x, t)H(x, y, t)\ , \quad (x, t)\in M^n\times (0, \infty)
}
\end{equation}
where $y\in M^n$ is fixed, $\Delta_x$ is the Laplacian operator with respect to $x$ and $g(t)$, and $R(x, t)$ is the scalar curvature at $x\in \big(M^n, g(t)\big)$.
}
\end{definition}

For the conjugate heat kernel $H(x, y, t)$ of backward Ricci flow, Perelman's \textbf{$\mathscr{W}$-entropy} is defined as 
\begin{align}
\mathscr{W}(g, f(x, y, t), t)= \int_{M^n} \Big[t\big(|\nabla f(x, y, t)|^2+ R(x, t)\big)+ f(x, y, t)- n \Big] H(x, y, t) d\mu_{g(t)}(x) \label{W}
\end{align}
where $H= (4\pi t)^{-\frac{n}{2}} e^{-f}$.

And Nash entropy for backward Ricci flow can be defined by imitating the linear heat equation case in \cite{Ni2}, 
\begin{align}
N(g, H, t)\vcentcolon = \Big(\int_{M^n} fH\Big)(t)- \frac{n}{2}\ , \label{def 3.1}
\end{align}

The limit of the reduced volume for backward Ricci flow is defined as 
\begin{align}
\widehat{V}_{\infty}(y, 0)\vcentcolon= \lim_{t\rightarrow \infty}  \int_{M^n} (4\pi t)^{-\frac{n}{2}} e^{-\ell (q, t)} d\mu_{g(t)}(q) \label{def V}
\end{align}
where $\ell(q, t)$ is the reduced distance with base-point $(y, 0)$ (check section $7$ of \cite{Pere} for the definition of the reduced distance).

Due to the surprising similarity between the entropy formula for the Ricci flow (\ref{W}), (\ref{def 3.1}), (\ref{def V}) and the entropy formula for the linear heat equation (\ref{1.1.NW}) (\ref{1.1.2}), the following question was asked by L. Ni in \cite{Ni2}:
\begin{question}\label{ques 1.1}
{Let $(M^n, g(t))_{t\in[0, \infty)}$ be a non-flat backward Ricci flow $\kappa$-solution with bounded nonnegative curvature operator $Rm$. One may ask if 
\begin{align}
\lim_{t\rightarrow \infty}\mathscr{W}(g, f(x, y, t), t)= \lim_{t\rightarrow \infty} N(g, H(x, y, t), t)= \ln \widehat{V}_{\infty}(y, 0) ? \label{1.2}
\end{align}
}
\end{question}

If we try to prove (\ref{1.2}) by imitating Ni's method mentioned above, we need to get a sharp bound of the conjugate heat kernel in the Ricci flow case. But it is much harder to get the sharp Gaussian bounds than the uniform Gaussian bounds for the conjugate heat kernel in the Ricci flow case, and we do not know the sharp bounds so far. 

On the other hand, by modifying the argument of Cao and Zhang in \cite{CZ}, we succeed in getting the uniform Gaussian bounds of the conjugate heat kernel $H(x, y, t)$ for Type I $\kappa$-solutions, where the coefficients of the bounds do not depend on time $t$, and the distance in exponential term and the volume term are with respect to the changing metric $g(t)$. Let us recall the definition of a Type I $\kappa$-solution,
\begin{definition}\label{def 1.3}
{A $\kappa$-solution on $[0, \infty)$ is called Type I if there exists a positive constant $C_1$ such that $|Rm(x, t)|\leq \frac{C_1}{1+ t}$ for any $t\in [0, \infty)$.
}
\end{definition} 

Also note that Cao and Zhang \cite{CZ} showed that the asymptotic limit of Type I $\kappa$-solution is a shrinking soliton. When $(M^n, g(t))_{t\in (0, \infty)}$ are shrinking soliton solutions to backward Ricci flow and $g(t)$ converges in the Gromov-Hausdorff sense as $t\searrow 0$ a metric cone $C$ which is smooth except at the vertex,  the above equation (\ref{1.2}) was proved by Cao, Hamilton and Ilmanen (see Section $3$ of \cite{CHI}). In the Ricci flow case, these provide us the analogue of Cheeger, Colding and Ding's results about heat equations on static manifolds.

Then following the similar strategy as the linear heat equation case, we get our main theorem as the following:
\begin{theorem}\label{thm 1.4}
{Let $(M^n, g(t))_{t\in [0, \infty)}$ be a non-flat Type I $\kappa$-solution to the backward Ricci flow for some $\kappa> 0$, and $Rm(x, t)\geq 0$ for all $(x, t)\in M^n\times[0, \infty)$. Then for $y\in M^n$, 
\begin{align}
\lim_{t\rightarrow \infty}\mathscr{W}(g, f(x, y, t), t)= \lim_{t\rightarrow \infty}N(g, H(x, y, t), t)= \ln \widehat{V}_{\infty}(y, 0) \nonumber
\end{align}
}
\end{theorem}
Because shrinking soliton solutions are obviously Type I ancient solutions, our theorem can be thought as a kind of generalization of the result of Cao, Hamilton and Ilmanen in Section $3$ of \cite{CHI}.

The paper is organized as follows: In Section \ref{Section 2}, we give the alternative proof of (\ref{1.1}). We prove the uniform (but not sharp) Gaussian bounds for conjugate heat kernel in Section \ref{section 3}, which is crucial for the later results. In Section \ref{section 4}, by similar argument as in Section $3$ of \cite{CHI}, (\ref{1.2}) is firstly proved for shrinking soliton solutions. Then we use the uniform Gaussian bound got in Section \ref{section 3} to reduce the proof of Theorem \ref{thm 1.4} about ancient solutions to shrinking soliton solutions.

\section{$\mathscr{W}$-entropy and reduced volume for the linear heat equation}\label{Section 2}

In this section, $(M^n, g)$ is a complete Riemannian manifold with $Rc\geq 0$ and maximum volume growth, namely $\theta_{\infty}= \lim_{r\rightarrow \infty} \frac{V_x(r)}{\omega_n r^n}> 0$.

For any increasing sequence $\{t_i\}$ with $\lim_{i\rightarrow \infty} t_i= \infty$, from Gromov's compactness theorem (see \cite{Gromov}), there exists a subsequence, also denoted as $\{t_i\}$, such that 
\begin{align}
(M^n, g_i, y)\rightarrow (M_\infty, y) \label{A.1}
\end{align}
where $g_i= t_i^{-1}g$, $M_{\infty}$ is some length space with measure $d\mu_{\infty}$ and the convergence is in the pointed Gromov-Hausdorff sense. From Theorem $7.6$ of \cite{CC}, $M_{\infty}$ is a metric cone, denoted as $C(X)$, where $X$ is a compact length space with measure $d\mu$. 

Define $H_i\doteqdot H_i(x, y, s)\doteqdot t_i^{\frac{n}{2}} H(x, y, t_i s)$, note that $H_i$ is a positive fundamental solution of heat equation on $(M^n, g_i)$. From Theorem $5.54$ and Theorem $6.20$ of \cite{Ding}, we get that for any $t> 0$,
\begin{align}
\lim_{i\rightarrow \infty}H_i(x, y, t)= H_{\infty}(x, y, t) \label{A.2}
\end{align}
the convergence is uniform in $C^0$-topology and also in $L^1$. And from $(6.23)$ of \cite{Ding}, we have
\begin{align}
H_{\infty}(x, y, t)= (4\pi t)^{-\frac{n}{2}} e^{-f_{\infty}}\ , \quad f_{\infty}(x, y, t)= \frac{r^2}{4t}+ \ln\Big(\frac{V(X)}{n\omega_n}\Big) \label{A.4}
\end{align}
where $r= d(x, y)$, $V(X)$ is the volume of $X$, and from the results in \cite{Colding0} it is easy to get 
\begin{align}
\frac{V(X)}{n\omega_n}= \frac{V_{M_{\infty}}\Big(B_{\infty}(y, 1)\Big)}{\omega_n}= \lim_{i\rightarrow \infty} \frac{V_{M^n}\Big(B(y, t_i)\Big)}{\omega_n t_i^n}= \theta_{\infty} \nonumber
\end{align}
From (\ref{A.4}) and the above, 
\begin{align}
f_{\infty}= \frac{r^2}{4t}+ \ln\theta_{\infty} \label{A.5}
\end{align}

On $M_{\infty}= C(X)$, fro the above explicit formulas we have the following lemma.
\begin{lemma}\label{lem A.1}
{\begin{align}
\int_{M_{\infty}} (f_{\infty}H_{\infty})(x, y, 1) d\mu_{\infty}(x)= \frac{n}{2}+ \ln\theta_{\infty} \nonumber
\end{align}
}
\end{lemma}

\pf
{From (\ref{A.5}), we get
\begin{align}
\Big(\int_{M_{\infty}} f_{\infty}H_{\infty} \Big)(1)= \int_{C(X)} \Big[\frac{r^2}{4}+ \ln\theta_{\infty}\Big]H_{\infty}(1)= \ln\theta_{\infty}+ \int_{C(X)}\frac{r^2}{4}\cdot H_{\infty}(1) \nonumber
\end{align}
On the other side, from (\ref{A.4}) we obtain
\begin{align}
\int_{C(X)} \frac{r^2}{4}\cdot H_{\infty}&= \int_0^{\infty}\int_X \Big[\frac{r^2}{4}\cdot (4\pi )^{-\frac{n}{2}} \exp{\Big(-\frac{r^2}{4}\Big)} \cdot \Big(\frac{n\omega_n}{V(X)}\Big)\Big] d\mu dr \nonumber \\
&= (n\omega_n \pi)\cdot (4\pi)^{-\frac{n}{2}- 1} \int_0^{\infty} r^{n+ 1} \exp{\Big(-\frac{r^2}{4}\Big)} dr \nonumber \\
&= \frac{n}{2}n\omega_n \pi^{-\frac{n}{2}} \Gamma (\frac{n}{2}+ 1)= \frac{n}{2} \nonumber
\end{align}
where $\Gamma$ is the Gamma function.

By all the above, the conclusion is proved.
}
\qed

The following result was firstly proved by L. Ni in \cite{Ni2}, \cite{Ni3} by different method.
\begin{prop}\label{prop A.2}
{\begin{align}
\lim_{t\rightarrow \infty} \mathscr{W}(f, t)= \lim_{t\rightarrow \infty} N(H, t)= \ln\theta_{\infty} \nonumber
\end{align}
}
\end{prop}

\pf
{We firstly show that $\lim_{t\rightarrow \infty} N(H, t)= \ln\theta_{\infty}$. From \cite{LY}, 
\begin{align}
\frac{C^{-1}(n)}{V_{y}(\sqrt{t})}\exp{\Big(-\frac{d^2(x, y)}{3t}\Big)} \leq H(x, y, t)\leq \frac{C(n)}{V_{y}(\sqrt{t})}\exp{\Big(-\frac{d^2(x, y)}{5t}\Big)} \label{A.2.1}
\end{align}
Recall $\lim_{r\rightarrow \infty} \frac{V_{x}(r)}{\omega_n r^n}= \theta_{\infty}> 0$ and the volume comparison theorem, from (\ref{A.2.1}) we have
\begin{align}
C^{-1}(n)\exp{\Big(-\frac{d^2(x, y)}{3t}\Big)}\leq t^{\frac{n}{2}} H\leq C(n)\exp{\Big(-\frac{d^2(x, y)}{5t}\Big)} \label{A.2.2}
\end{align}

Hence for any $b\geq 1$,
\begin{align}
\Big|\int_{M^n\backslash B(y, b\sqrt{t})} (fH)(x, y, t) d\mu_{g}(x)\Big|& \leq \int_{M^n\backslash B(y, b\sqrt{t})} \Big|\ln \big[(4\pi t)^{\frac{n}{2}} H\big]\Big|\cdot H \nonumber \\
&\leq \int_{M^n\backslash B(y, b\sqrt{t})} \Big(C+ \frac{d^2(x, y)}{3t}\Big)\cdot H(x, y, t) d\mu_{g}(x) \label{A.2.3}
\end{align}
where $C= C(n)$, and in the last inequality we used (\ref{A.2.2}). In the rest of the proof $C= C(n)$ if not specifically mentioned,

Using (\ref{A.2.1}) and do integration by parts, we obtain that
\begin{align}
&\int_{M^n\backslash B(y, b\sqrt{t})} \Big(C+ \frac{d^2(x, y)}{3t}\Big)\cdot H(x, y, t) d\mu_{g}(x) \nonumber \\
& \leq \frac{C}{V_{y}(\sqrt{t})} \int_{b\sqrt{t}}^{\infty} \exp{\Big(-\frac{r^2}{5t}\Big)}\Big(C+ \frac{r^2}{t}\Big) dV_{y}(r) \nonumber \\
& \leq C\int_{b\sqrt{t}}^{\infty} \frac{V_y(r)}{V_y(\sqrt{t})}  \exp{\Big(-\frac{r^2}{5t}\Big)}\Big(\frac{C}{5}+ \frac{r^2}{5t}\Big) d\Big(\frac{r^2}{t}\Big) \nonumber \\
& \leq C\int_{b^2}^{\infty} s^{\frac{n}{2}} \Big(\frac{s}{5}+ C\Big) e^{-\frac{s}{5}} ds \label{A.2.4}
\end{align}
in the last equality the volume comparison theorem is used.

It is easy to get 
\begin{align}
\int_{b^2}^{\infty} s^{\frac{n}{2}}\Big(\frac{s}{5}+ C\Big) e^{-\frac{s}{5}} ds\leq C\int_{b^2}^{\infty} e^{-\frac{s}{10}} ds\leq C\cdot e^{-\frac{b^2}{10}} \label{A.2.5}
\end{align}

From (\ref{A.2.3}), (\ref{A.2.4}) and (\ref{A.2.5}), it follows that
\begin{align}
\Big|\int_{M^n\backslash B(y, b\sqrt{t})} (fH)(x, y, t) d\mu_{g}(x)\Big| \leq C\cdot e^{-\frac{1}{10}b^2} \label{A.2.6}
\end{align}

Then 
\begin{align}
\lim_{i\rightarrow \infty} \Big(\int_{M^n} fH\Big)(t_i)&= \lim_{i\rightarrow \infty} \Big(\int_{M^n\backslash B(y, b\sqrt{t_i})} fH\Big)(t_i)+ \lim_{i\rightarrow \infty} \Big(\int_{B_{g_i}(y, b)} f_i H_i \Big)(1) \nonumber \\ 
&\geq -C\cdot e^{-\frac{1}{10}b^2} + \Big(\int_{B_{\infty}(y, b)} f_{\infty} H_{\infty} \Big)(1) \nonumber 
\end{align}
where $f_i$ is defined by $H_i (s)= (4\pi s)^{-\frac{n}{2}} \exp{\Big(-f_{i}\Big)}$, and in the last inequality we used (\ref{A.2.6}) and (\ref{A.2}). Let $b\rightarrow \infty$ in the above, we deduce
\begin{align}
\lim_{i\rightarrow \infty}\Big(\int_{M^n} fH\Big)(t_i)\geq \Big(\int_{M_{\infty}} f_{\infty} H_{\infty}\Big) (1) \nonumber
\end{align}

On the other hand, similarly using (\ref{A.2.6}) and (\ref{A.2}), we arrive at
\begin{align}
\lim_{i\rightarrow \infty}\Big(\int_{M^n} fH\Big)(t_i)\leq \Big(\int_{M_{\infty}} f_{\infty} H_{\infty}\Big) (1) \nonumber
\end{align}

By all the above and Lemma \ref{lem A.1}, 
\begin{align}
\lim_{i\rightarrow \infty} \Big(\int_{M^n} fH\Big)(t_i)= \Big(\int_{M_{\infty}} f_{\infty} H_{\infty}\Big) (1)= \frac{n}{2}+ \ln\theta_{\infty} \nonumber 
\end{align}

We know that $\theta_{\infty}$ is independent of the choice of the sequence $\{t_i\}$, hence 
\begin{align}
\lim_{t\rightarrow \infty} \Big(\int_{M^n} fH\Big)(t)= \frac{n}{2}+ \ln\theta_{\infty} \nonumber
\end{align}
and it is equivalent to
\begin{align}
\lim_{t\rightarrow \infty} N(H, t)= \ln\theta_{\infty} \label{A.2.7} 
\end{align}

From (\ref{A.2.7}), we can get that $\Big|N(H, 2t)- N(H, t)\Big|\leq \epsilon$ for $t\gg 1$. This implies that there exists $\{t_i\}$ such that $t_i\frac{\partial}{\partial t}N(H, t_i)\rightarrow 0$ as $t_i\rightarrow \infty$. From the monotonicity of $\mathscr{W}(f, t)$ shown in \cite{Ni1} and (\ref{A.2.7}), we finally get
\begin{align}
\lim_{t\rightarrow \infty}\mathscr{W}(f, t)= \lim_{i\rightarrow \infty}\mathscr{W}(f, t_i)= \lim_{i\rightarrow \infty}\Big[t_i\frac{\partial}{\partial t}N(H, t_i)+ N(H, t_i)\Big]= \ln\theta_{\infty} \nonumber
\end{align}

The conclusion is proved.
}
\qed

\section{Gaussian bounds of the conjugate heat kernel in the Ricci flow} \label{section 3}

In this section, we assume that $(M^n, g(t))_{t\in [0, \infty)}$ is a non-flat Type I $\kappa$-solution to the backward Ricci flow for some $\kappa> 0$ and 
$|Rm(x, t)|\leq \frac{C_1}{1+ t}$, and $H(x, y, t)$ is the conjugate heat kernel defined in the definition \ref{def chk}.

We will use the following result due to Cao and Zhang repeatedly.
\begin{lemma}[Lemma $4.1$ of \cite{CZ}]\label{lem 4.1 of CZ}
{\begin{equation}\nonumber
{C^{-1}t^{-\frac{n}{2}} \leq H(x, y, t)\leq Ct^{-\frac{n}{2}}
}
\end{equation}
where $C= C(C_1, n, \kappa)$ and $(x, y, t)\in M^n\times M^n\times (0, \infty)$.
}
\end{lemma}

We have the following curvature derivative estimates from the assumption and Shi's estimates.
\begin{lemma}\label{lem 2.1}
{\begin{equation}\label{2.1}
{|\nabla Rm(x, t)|\leq C t^{-\frac{3}{2}}\ , \quad for \ \  (x, t)\in M^n\times (0, \infty)
}
\end{equation}
where $C= C(C_1, n)$.
}
\end{lemma}

\pf
{Fix $T> 0$, define $\tilde{g}(t)= g(2T- t)$, then 
\begin{equation}\nonumber
{|Rm|_{\tilde{g}(t)}\leq \frac{C_1}{1+ (2T- t)}\leq \frac{C_1}{T}\ , \quad for \ \ t\in[0, T]
}
\end{equation}

We use Shi's global derivative estimates Theorem $14.5$ in \cite{RFAA}, choose $\alpha= C_1$, $K= \frac{C_1}{T}$ there, we have
\begin{equation}\nonumber
{|Rm|_{\tilde{g}(t)}\leq K\ , \quad for \ \ t\in[0, \frac{\alpha}{K}]
}
\end{equation}
hence there exists $C= C(n, \max\{\alpha, 1\})= C(C_1, n)$, such that 
\begin{equation}\nonumber
{|\nabla Rm|_{\tilde{g}(t)}\leq \frac{C}{\sqrt{t}}K\ , \quad for \ \ t\in(0, \frac{\alpha}{K}]= (0, T]
}
\end{equation}

Let $t= T$, we get
\begin{equation}\nonumber
{|\nabla Rm|_{g(T)}= |\nabla Rm|_{\tilde{g}(T)}\leq \frac{C}{T^{\frac{3}{2}}}
}
\end{equation}

Because $T> 0$ is chosen freely, the lemma is proved.
}
\qed

We can use the curvature estimates obtained and the result in \cite{RFAA} to get the following gradient estimate for $H(x, y, t)$, which is different from the gradient estimate proved in \cite{Zhang}.
\begin{prop}\label{prop 2.2}
{\begin{equation}\label{2.3}
{\Big|\frac{\nabla H}{H}\Big|(x, y, t)\leq \frac{C}{\sqrt{t}}
}
\end{equation}
where $\nabla$ is with respect to $x$, and $C= C(C_1, n, \kappa)$.
}
\end{prop}

\pf
{Fix $T> 0$, define $\tilde{g}(t)= g(\frac{T}{2}+ t)$, use Lemma \ref{lem 2.1} we get
\begin{equation}\nonumber
{|\nabla Rm|_{\tilde{g}(t)}\leq \frac{C}{(t+ \frac{T}{2})^{\frac{3}{2}}}\leq \Big(\frac{C}{T}\Big)^{\frac{3}{2}}\ , \quad t\in [0, \frac{T}{2}]
}
\end{equation}

Define $\tilde{H}(x, y, t)= H(x, y, t+ \frac{T}{2})$, then $\tilde{H}(t)$ is the solution to (\ref{2.2}) with respect to metric $\tilde{g}(t)$. From Lemma \ref{lem 4.1 of CZ},
\begin{equation}\nonumber
{C^{-1}T^{-\frac{n}{2}}\leq \tilde{H}(x, y, t)\leq C T^{-\frac{n}{2}}\ , \quad t\in [0, \frac{T}{2}]
}
\end{equation}
where $C= C(n, C_1, \kappa)$.

Using $Rc\geq -\frac{C}{1+ t}$ and Theorem $16.52$ in \cite{RFAA} for $(M^n, \tilde{g}(t))$, $t\in [0, \frac{T}{2}]$, we get 
\begin{equation}\nonumber
{\Big|\frac{\nabla \tilde{H}}{\tilde{H}}\Big|^2(x, y, t)\leq \Big[1+ \ln \big(\frac{CT^{-\frac{n}{2}}}{\tilde{H}}\big)\Big]^2\Big(\frac{C}{t}+ 2\frac{C}{T}\Big)
}
\end{equation}
where $\nabla$ is with respect to $x$. Let $t= \frac{T}{2}$ in the above inequality, we get
\begin{equation}\nonumber
{\Big|\frac{\nabla H}{H}\Big|^2(x, y, T)= \Big|\frac{\nabla \tilde{H}}{\tilde{H}}\Big|^2(x, y, \frac{T}{2})\leq \frac{C}{T}
}
\end{equation}
where $C= C(C_1, n, \kappa)$. Because $T> 0$ is chosen freely, the conclusion is proved.
}
\qed

A corollary is the following Harnack inequality.
\begin{cor}\label{cor 2.3}
{\begin{equation}\label{2.4}
{H(x, y, t)\leq H(z, y, t)\exp{\Big(\frac{Cd(x, z, t)}{\sqrt{t}}\Big)}
}
\end{equation}
where $C= C(C_1, n , \kappa)$.
}
\end{cor}

\pf
{Fix $x$ and $z$, then integrate (\ref{2.3}) along a minimal geodesic connecting $x$ and $z$.
}
\qed

Now we are ready to prove the main theorem of this section.
\begin{theorem}\label{thm 2.4}
{There exist positive constants $\Lambda_i= \Lambda_i(C_1, n, \kappa)$, $i= 1, 2, 3, 4$, such that
\begin{equation}\label{2.5}
{\frac{\Lambda_1}{V_{g(t)}\Big(B(x, \sqrt{t}, t)\Big)} \exp{\Big(-\frac{\Lambda_2 d^2}{t}\Big)}\leq H\leq \frac{\Lambda_3}{V_{g(t)}\Big(B(x, \sqrt{t}, t)\Big)} \exp{\Big(-\frac{\Lambda_4 d^2}{t}\Big)}
}
\end{equation}
where $B(x, r, t)= \{z|\ d_{g(t)}(z, x)\leq r\}$ for any $r> 0$, $H= H(x, y, t)$, $d= d(x, y, t)$ and $(x, y, t)\in M^n\times M^n\times (0, \infty)$.
}
\end{theorem}

\begin{remark}\label{rem 2.5}
{We use the method of Grigor'yan (see \cite{Grig}) as the proof of Theorem $3.1$ in \cite{CZ}, but we estimate the conjugate heat kernel directly, which is different from Cao and Zhang's strategy in \cite{CZ}. In \cite{Ni-add}, the related two sided bounds of heat kernel was got in Theorem $1.3$ there. The key difference between Theorem $1.3$ in \cite{Ni-add} and Theorem \ref{thm 2.4} above, is that the volume terms $V_{g(t)}$ and distance terms $d$ appearing in our two sided bounds of $H$ depend on the evolving metric $g(t)$ (not the fixed metric $g(0)$), which is crucial for the proof of Theorem \ref{thm 1.4} in Section \ref{section 4}.
}
\end{remark}

\pf
{$\mathbf{Step\ 1}$. Pick a weight function $e^{\xi(x, t)}$ which will be specified later. Using integration by parts, we can get 
\begin{align}
\frac{\partial}{\partial t}\int_{M^n} H^2(x, y, t)e^{\xi(x, t)} d\mu_{g(t)}(x)&= \int_{M^n} H^2 e^{\xi}\xi_t+ (2H\Delta H- RH^2)e^{\xi} \nonumber \\
&\leq \int_{M^n} \Big(\xi_t+ \frac{1}{2}|\nabla \xi|^2- R\Big)H^2 e^{\xi} \label{2.4.1}
\end{align}
Because $M^n$ can be non-compact, one needs to justify integration by parts near infinity in (\ref{2.4.1}). For fixed $t$, $H(x, y, t)$ has a generic Gaussian upper bound with coefficients depending on $t$, curvature tensor and their derivatives, as shown in \cite{CTY}. Since the curvatures are all bounded, from Proposition \ref{prop 2.2} and volume comparison theorem, the term $\int_{\partial B(r)}|\nabla H|H e^{\xi}\rightarrow 0$ as $r\rightarrow \infty$, this justifies the integration by parts in (\ref{2.4.1}).

We choose 
\begin{equation}\nonumber
\xi(x, t)= \left\{
\begin{array}{rl}
-\frac{(\iota- d(x, y, t))^2}{4(s_0- t)}\ , \quad &d(x, y, t)\leq \iota\ ;\\
0\ , \quad &d(x, y, t)> \iota \ .
\end{array} \right.
\end{equation}
where $\iota> 0$ and $s_0> t> 0$.

Then for $x\in B(y, \iota, t)$, we have 
\begin{equation}\nonumber
{\xi_t+ \frac{1}{2}|\nabla \xi|^2= -\frac{1}{8}\cdot \Big(\frac{\iota- d}{s_0- t}\Big)^2+ \frac{1}{2}\cdot \frac{\iota- d}{s_0- t}\cdot \frac{\partial}{\partial t}d(x, y, t)
}
\end{equation}

By Lemma $8.3$ $(b)$ of \cite{Pere} and $|Rc|(z, t)\leq \frac{C}{1+ t}$ for any $z\in M^n$, where $C= C(C_1, n)$, we get 
\begin{equation}\nonumber
{\frac{\partial}{\partial t}d(x, y, t)\leq 2(n- 1)\Big(\frac{2}{3}\frac{C}{1+ t}r_0+ r_0^{-1}\Big)\ , \quad for \ any \ r_0> 0
}
\end{equation}
we can choose $r_0= \sqrt{1+ t}$, then $\frac{\partial}{\partial t}d\leq \frac{C}{\sqrt{1+ t}}$, where $C= C(C_1, n)$. It follows that 
\begin{align}
\xi_t+ \frac{1}{2}|\nabla \xi|^2&\leq -\frac{1}{8}\cdot \Big(\frac{\iota- d}{s_0- t}\Big)^2+ \frac{1}{2}\cdot \frac{\iota- d}{s_0- t}\cdot \frac{C}{\sqrt{1+ t}} \nonumber \\
&= -\frac{1}{8}\Big(\frac{\iota- d}{s_0- t}- 2\frac{C}{\sqrt{1+ t}}\Big)^2+ \frac{1}{2}\frac{C^2}{1+ t}\leq \frac{C}{1+ t} \nonumber
\end{align}

Combining with $Rc\geq -\frac{C}{1+ t}$, from (\ref{2.4.1}) we get 
\begin{equation}\nonumber
{\frac{\partial}{\partial t}\Big(\int_{M^n} H^2 e^{\xi}\Big)\leq \frac{C_2}{1+ t} \Big(\int_{M^n} H^2 e^{\xi}\Big)
}
\end{equation}
where $C_2= C_2(C_1, n)$. Hence
\begin{equation}\label{2.4.2}
{\Big(\int_{M^n} H^2 e^{\xi}\Big)(s_1)\leq \Big(\int_{M^n} H^2 e^{\xi}\Big)(s_2)\Big(\frac{1+ s_1}{1+ s_2}\Big)^{C_2}
}
\end{equation}
where $s_0$, $s_1$, $s_2$ are any positive constants satisfying $s_0> s_1> s_2> 0$.

$\mathbf{Step\ 2}$. Define 
\begin{equation}\nonumber
{I_{\iota}(t)= \int_{M^n\backslash B(y, \iota, t)} H^2(x, y, t) d\mu_{g(t)}(x)
}
\end{equation}

Choose any $0< \rho< \iota$, using (\ref{2.4.2}), we have 
\begin{align}
I_{\iota}(s_1)&= \int_{M^n\backslash B(y, \iota, t)} H^2(x, y, s_1) \leq \int_{M^n} H^2(x, y, s_1)e^{\xi(x, s_1)} \nonumber \\
&\leq \Big[\int_{M^n} H^2 e^{\xi}(s_2)\Big] \Big(\frac{1+ s_1}{1+ s_2}\Big)^{C_2} \nonumber \\
&\leq \Big[I_{\rho}(s_2)+ \exp{\Big(-\frac{(\iota- \rho)^2}{4(s_0- s_2)}\Big)}\int_{B(y, \rho , s_2)} H^2(s_2) \Big] \Big(\frac{1+ s_1}{1+ s_2}\Big)^{C_2} \nonumber \\
&\leq \Big[I_{\rho}(s_2)+ C_3 s_2^{-\frac{n}{2}}\exp{\Big(-\frac{(\iota- \rho)^2}{4(s_0- s_2)}\Big)} \Big] \Big(\frac{1+ s_1}{1+ s_2}\Big)^{C_2} \nonumber 
\end{align}
where $C_3= C(C_1, n, \kappa)$, and in the last inequality we used Lemma \ref{lem 4.1 of CZ} and $\int_{M^n} H\equiv 1$.

Let $s_0\rightarrow s_1$ in the above inequality,
\begin{align}
I_{\iota}(s_1)\leq \Big[I_{\rho}(s_2)+ C_3 s_2^{-\frac{n}{2}}\exp{\Big(-\frac{(\iota- \rho)^2}{4(s_1- s_2)}\Big)} \Big] \Big(\frac{1+ s_1}{1+ s_2}\Big)^{C_2} \label{2.4.3}
\end{align}
Note (\ref{2.4.3}) holds for any $s_1> s_2> 0$, $\iota> \rho> 0$. 

Now we define 
\begin{equation}\nonumber
{t_k= ta^{- k}\ , \quad  r_k=\Big(\frac{1}{2}+ \frac{1}{k+ 2}\Big)r\ , \quad k= 0, 1, 2, \cdots
}
\end{equation}
where $a> 1$ is a constant to be chosen later. Let $s_1= t_k$, $s_2= t_{k+ 1}$, $\iota= r_k$ and $\rho= r_{k+ 1}$, applying (\ref{2.4.3}) we can get
\begin{align}
I_{r_k}(t_k)\leq \Big[I_{r_k+ 1}(t_{k+ 1})+ C_3 t_{k+ 1}^{- \frac{n}{2}} \exp{\Big(-\frac{(r_k- r_{k+ 1})^2}{4(t_k- t_{k+ 1})}\Big)} \Big]\cdot \Big(\frac{1+ t_{k}}{1+ t_{k+ 1}}\Big)^{C_2} \label{2.4.4}
\end{align} 

After applying iteration to (\ref{2.4.4}), we obtain 
\begin{align}
I_{r}(t)= I_{t_0}(t_0)\leq \Big(\frac{1+ t_0}{1+ t_{k+ 1}}\Big)^{C_2}\cdot I_{r_{k+ 1}}(t_{k+ 1})+ C_{3} t^{-\frac{n}{2}} \sum_{j= 1}^{k+ 1} a^{\big(\frac{j}{2}+ C_2\big)n} \exp{\Big(-\frac{r^2a^j}{4(j+ 3)^4(a- 1)t}\Big)} \label{2.4.5}
\end{align} 

When $k\rightarrow \infty$, $t_k\rightarrow 0$ and $H(x, y, t_k)\rightarrow \delta_y(x)$ which is concentrated at the point $y$. Hence $\lim_{k\rightarrow \infty}I_{r_k}(t_k)= 0$. Let $k\rightarrow \infty$ in (\ref{2.4.5}), we get
\begin{align}
I_r(t)\leq C_3 t^{-\frac{n}{2}} \sum_{j= 1}^{\infty} a^{\big(\frac{j}{2}+ C_2\big)n} \exp{\Big(-\frac{r^2a^j}{4(j+ 3)^4(a- 1)t}\Big)} \nonumber
\end{align}
By taking $r^2\geq \frac{1}{4}t$ and making the constant $a$ sufficiently large, it leads to 
\begin{align}
\int_{M^n\backslash B(y, r, t)} H^2(t)= I_r(t)\leq C_{4}t^{-\frac{n}{2}} \exp{\Big(-\frac{C_5r^2}{t}\Big)} \nonumber
\end{align}
where $C_4= C_4(C_1, n , \kappa)$ and $C_5= C_5(C_1, n)$. 

Using Lemma \ref{lem 4.1 of CZ}, we can get
\begin{align}
\int_{M^n\backslash B(y, r, t)} H(x, y, t)\leq C_{6} \exp{\Big(-\frac{C_5r^2}{t}\Big)} \label{2.4.7}
\end{align}
where $C_6= C_6(C_1, n, \kappa)$ and $r\geq \frac{1}{2}\sqrt{t}$.

$\mathbf{Step\ 3}$. Let $x_0\in M^n$, there are two cases.

If $d(x_0, y, t)\geq \sqrt{t}$, then $B(x_0, \frac{1}{2}\sqrt{t}, t)\subset M^n\backslash B(y, r, t)$, where $r= \frac{1}{2}d(x_0, y, t)$. Then from (\ref{2.4.7}), there exists $z_0\in B(x_0, \frac{1}{2}\sqrt{t}, t)$ such that 
\begin{align}
H(z_0, y, t)V_{g(t)}\Big(B(x_0, \frac{1}{2}\sqrt{t}, t)\Big)\leq C_6 \exp{\Big(-\frac{C_5r^2}{t}\Big)} \nonumber
\end{align} 

Note the following fact
\begin{align}
\frac{V_{K}(r_1)}{V_{K}(r_2)}\leq \Big(\frac{r_1}{r_2}\Big)^{n}e^{n\sqrt{K}r_1} \label{fact}
\end{align}
where $V_{K}(r_i)$ denotes the volume of a ball of radius $r_i$ in the constant curvature $-K$ $n$-dimensional space form, and $K> 0$.

From the above fact, by the classical volume comparison theorem and $Rc\geq -\frac{C}{1+ t}$, these imply that 
\begin{align}
H(z_0, y, t)\leq C_7\Big[V_{g(t)}\Big(B(x_0, \sqrt{t}, t)\Big)\Big]^{-1} \exp{\Big(-\frac{C_8d^2(x_0, y, t)}{t}\Big)}
\end{align}
where $C_7= C_7(C_1, n, \kappa)$ and $C_8= C_8(C_1, n)$.

By Corollary \ref{cor 2.3}, using $d(x_0, z_0, t)\leq \frac{1}{2}\sqrt{t}$ we have
\begin{align}
H(x_0, y, t)&\leq H(z_0, y, t)\cdot \exp{\Big(\frac{C_9 d(x_0, z_0, t)}{\sqrt{t}}\Big)} \nonumber \\
&\leq C_{10}\Big[V_{g(t)}\Big(B(x_0, \sqrt{t}, t)\Big)\Big]^{-1} \exp{\Big(-\frac{C_8d^2(x_0, y, t)}{t}\Big)} \label{2.4.8}
\end{align}
where $C_{10}= C_{10}(C_1, n , \kappa)$.

If $d(x_0, y, t)\leq \sqrt{t}$, by Lemma \ref{lem 4.1 of CZ}, using the volume comparison theorem and $Rc\geq -\frac{C}{1+ t}$ again, there exists $C_{11}= C_{11}(C_{10}, C_8, n)= C_{11}(C_1, n, \kappa)$ such that
\begin{align}
H(x_0, y , t)\leq Ct^{-\frac{n}{2}}\leq  C_{11} \Big[V_{g(t)}\Big(B(x_0, \sqrt{t}, t)\Big)\Big]^{-1} \exp{\Big(-\frac{C_8 d^2(x_0, y, t)}{t}\Big)} \label{2.4.9}
\end{align} 

Define $C_{12}= \max{\{C_{10}, C_{11}\}}$, then from (\ref{2.4.8}) and (\ref{2.4.9}) we get
\begin{align}
H(x_0, y , t)\leq Ct^{-\frac{n}{2}}\leq  C_{12} \Big[V_{g(t)}\Big(B(x_0, \sqrt{t}, t)\Big)\Big]^{-1} \exp{\Big(-\frac{C_8 d^2(x_0, y, t)}{t}\Big)} \nonumber
\end{align} 

Since $x_0$ is arbitrary, this proves the desired upper bound.

$\mathbf{Step\ 4}$. Next we show that a lower bound follows from the upper bound. From (\ref{2.4.7}), we get 
\begin{align}
\int_{B(y, b\sqrt{t}, t)} H \geq 1- C_6\exp{\Big(-C_5 b^2\Big)} \label{2.4.9.1}
\end{align}
where $b\geq 1$ is a constant to be chosen later. Hence there exists $x_1\in B(y, b\sqrt{t}, t)$ such that 
\begin{align}
H(x_1, y, t)\geq \Big[V_{g(t)}\Big(B(y, b\sqrt{t}, t)\Big)\Big]^{-1} \Big[1- C_6\exp{\Big(-C_5 b^2\Big)} \Big]\nonumber
\end{align}

For any $x_2\in M^n$, from Corollary \ref{cor 2.3},
\begin{align}
H(x_2, y, t)&\geq H(x_1, y, t)\exp{\Big(-\frac{C_9 d(x_1, x_2, t)}{\sqrt{t}}\Big)} \nonumber \\
&\geq \Big[V_{g(t)}\Big(B(y, b\sqrt{t}, t)\Big)\Big]^{-1} \Big[1- C_6\exp{\Big(-C_5 b^2\Big)} \Big] \exp{\Big(-\frac{C_9 d(x_1, x_2, t)}{\sqrt{t}}\Big)} \label{2.4.10}
\end{align}

Note 
\begin{align}
d^2(x_1, x_2, t)\leq \Big[d(x_1, y, t)+ d(x_2, y, t)\Big]^2\leq 2b^2t+ 2d^2(x_2, y, t) \nonumber
\end{align}
then
\begin{align}
\frac{d(x_1, x_2, t)}{\sqrt{t}} \leq \frac{d^2(x_1, x_2, t)}{t}+ 1\leq \frac{2d^2(x_2, y, t)}{t}+ (2b^2+ 1) \label{2.4.11}
\end{align}

From (\ref{fact}), combining with the volume comparison theorem and $Rc\geq -\frac{C}{1+ t}$, we can get 
\begin{align}
V_{g(t)}\Big(B(y, b\sqrt{t}, t)\Big)&\leq V_{g(t)}\Big(B(x_2, d(x_2, y, t)+ b\sqrt{t}, t)\Big) \nonumber \\
&\leq C_{13} \exp{\Big(\frac{C_{14}d(x_2, y, t)}{\sqrt{t}} \Big)} \cdot V_{g(t)}\Big(B(x_2, \sqrt{t}, t)\Big) \cdot \Big(\frac{d(x_2, y, t)+ b\sqrt{t}}{\sqrt{t}}\Big)^n \label{2.4.12}
\end{align}
where $C_{13}= C_{13}(C_1, n, b)$ and $C_{14}= C_{14}(C_1, n)$.

Choose the constant $b$ large enough such that $1- C_6\exp{\Big(-C_5 b^2\Big)}\geq \frac{1}{2}$, then from (\ref{2.4.10}), (\ref{2.4.11}) and (\ref{2.4.12}),
\begin{align}
H(x_2, y, t)&\geq C_{15} \Big[V_{g(t)}\Big(B(x_2, \sqrt{t}, t)\Big)\Big]^{-1}\cdot \Big(\frac{d(x_2, y, t)}{\sqrt{t}}+ b\Big)^{-n} \exp{\Big(-\frac{C_{16}d^2(x_2, y, t)}{t}\Big)} \nonumber \\
&\geq C_{17}\Big[V_{g(t)}\Big(B(x_2, \sqrt{t}, t)\Big)\Big]^{-1} \cdot \exp{\Big(-\frac{C_{18}d^2(x_2, y, t)}{t}\Big)} \nonumber
\end{align}

Since $x_2$ is arbitrary, this is a lower bound which matches the upper bound except for constant coefficients $\Lambda_i$.
}
\qed

\section{The limit of reduced volume and $\mathscr{W}$-entropy for the Ricci flow} \label{section 4}

Assume that $(M^n, g(t))_{t\in [0, \infty)}$ is a non-flat Type I $\kappa$-solution to the backward Ricci flow for some $\kappa> 0$, $H(x, y, t)$ is the fundamental solution to (\ref{2.2}) where $y\in M^n$ is a fixed point. By the $\kappa$-noncollapsed assumption and curvature bound $|Rm(\cdot, t)|\leq \frac{C_1}{1+ t}$, it follows from Hamilton's compactness theorem (see \cite{Ham}): for any increasing sequence $\{t_i\}$ with $\lim_{i\rightarrow \infty}t_i= \infty$, there exist a subsequence, also denoted as $\{t_i\}$, such that the following statement holds:

The pointed manifolds $(M^n, g_i(s), y)$ with metrics $g_i(s)\doteqdot t_i^{-1} g(t_i s)$ converges to a pointed manifold $(M_{\infty}, g_{\infty}(s), y)$ in $C^{\infty}_{loc}$-topology, where $s\in (0, \infty)$. 

It was shown in \cite{CZ} that $(M_{\infty}, g_{\infty}(s))$ is a gradient shrinking soliton. For completeness and later use, we give the details here following the argument in \cite{CZ}.

Define $H_i\doteqdot H_{i}(x, y, s)\doteqdot t_{i}^{\frac{n}{2}}H(x, y, t_is)$. By Lemma \ref{lem 4.1 of CZ}, for fixed $s> 0$, there exists a uniform positive constant $U_0$, such that $H_i(x, y, s)\leq U_0$ for all $i= 1, 2, \cdots$ and $x\in M^n$. Note that $H_i$ is a positive fundamental solution of the conjugate heat equation on $(M^n, g_i(s))$, i.e., 
\begin{align}
\frac{\partial}{\partial t}H_i= \Delta_{g_i}H_i- R_{g_i}H_i \nonumber
\end{align}

For any compact time interval in $(0, \infty)$, $H_i$ are uniformly bounded, moreover, $R_{g_i}$ and $Rm_{g_i}$ are uniformly bounded. It follows from the standard parabolic theory that $H_i$ is H\"older continuous uniformly with respect to $g_i$. Hence there exist a subsequence, still denoted as $\{H_{i}(x, y, s)\}$, which converges to a $C^{\alpha}_{loc}$-topology sense. 

It is easy to see that $H_{\infty}$ is a weak solution of the conjugate heat equation on $(M_{\infty}, g_{\infty}(s))$. By standard parabolic theory and the boundedness of $H_{\infty}$ on compact time interval, $H_{\infty}$ is a smooth solution of the conjugate heat equation on $(M_{\infty}, g_{\infty}(s))_{s\in (0, \infty)}$.

By Lemma \ref{lem 4.1 of CZ}, $H(x, y, t_is)\geq \frac{C}{(t_i s)^{\frac{n}{2}}}$, which is equivalent to $H_i(x, y, s)\geq Cs^{-\frac{n}{2}}$. Hence $H_{\infty}(s)\geq Cs^{-\frac{n}{2}}> 0$, it yields that $H_{\infty}$ is positive everywhere when $s> 0$.

We define 
\begin{align}
\mathscr{W}_i(s)\doteqdot \mathscr{W}(g_i, f_i, s)= \int_{M^n} \Big[s(|\nabla f_i|^2+ R_i)+ f_i- n\Big]H_i d\mu_{g_i(s)}(x) \nonumber
\end{align}
where $H_i(x, y, s)= (4\pi s)^{-\frac{n}{2}} e^{-f_i(x, y, s)}$, and $R_i$ is the scalar curvature with respect to $g_i$.

Because $M^n$ may be noncompact, one needs to justify that the integral $\mathscr{W}_i(s)$ is finite. This can be deduced from Lemma \ref{lem 4.1 of CZ}, Proposition \ref{prop 2.2} and $|Rm(\cdot, t)|\leq \frac{C_1}{1+ t}$ easily.

From Corollary $2.5$ in \cite{Chen}, $R\geq 0$ for any complete ancient solution to the Ricci flow, then for fixed $s> 0$, 
\begin{align}
\mathscr{W}_i(s)&\geq \Big(\int_{M^n} f_iH_i\Big)- n= -n -\int_{M^n} \ln \Big[ (4\pi)^{\frac{n}{2}} s^{\frac{n}{2}} H_i(s)\Big] H_i \nonumber \\
&\geq -n+ \int_{M^n} C\cdot H_i= C \label{3.1}
\end{align}
where $C$ is independent of $i$, in the last inequality we used Lemma \ref{lem 4.1 of CZ}.

Recall that $\mathscr{W}$ is invariant under proper scaling, 
\begin{align}
\mathscr{W}_i(s)= \mathscr{W}(g_i, f_i, s)= \mathscr{W}(g, f, t_i s) \nonumber
\end{align}

According \cite{Pere}, for fixed $s> 0$, $\mathscr{W}_i(s)= \mathscr{W}(g, f, t_i s)$ is a non-increasing sequence of $i$. By (\ref{3.1}), there exists a function $\mathscr{W}_{\infty}(s)$ such that 
\begin{align}
\mathscr{W}_{\infty}(s)\doteqdot \lim_{i\rightarrow \infty} \mathscr{W}_i(s)= \lim_{i\rightarrow \infty}\mathscr{W}(g, f, t_i s) \nonumber
\end{align}
Note that $\mathscr{W}_{\infty}(s)$ is independent of the choice of $\{t_i\}$ by the monotonicity of $\mathscr{W}(g, f, s)$.

For any fixed $s_0\in (0, \infty)$, we can find a subsequence $\{t_{m_i}\}_{i=1}^{\infty}$ tending to infinity such that 
\begin{align}
\mathscr{W}(g, f, t_{m_i} s_0)\geq \mathscr{W}(g, f, t_{m_i}(s_0+ 1))\geq \mathscr{W}(g, f, t_{m_{i+ 1}} s_0) \nonumber
\end{align}
Since $\lim_{i\rightarrow \infty}\mathscr{W}(g, f, t_{m_i} s_0)= \lim_{i\rightarrow \infty} \mathscr{W}(g, f, t_{m_{i+ 1}}s_0)= \mathscr{W}_{\infty} (s_0)$, we get 
\begin{align}
\lim_{i\rightarrow \infty} \Big[\mathscr{W}(g, f, t_{m_i}s_0)- \mathscr{W}(g, f, t_{m_i}(s_0+ 1))\Big]= 0 \nonumber
\end{align}
That is 
\begin{align}
\lim_{i\rightarrow \infty}\Big[\mathscr{W}_{m_i}(s_0)- \mathscr{W}_{m_i}(s_0+ 1)\Big]= 0 \label{3.2}
\end{align}

According \cite{Pere}, 
\begin{align}
\frac{d}{d s}\mathscr{W}_k(s)= -2\int_{M^n} \Big|Rc_{g_k}+ Hess_{g_k}f_k- \frac{1}{2s}g_k \Big|^2 H_k d\mu_{g_k}(s)  \nonumber
\end{align}
Integrate it from $s_0$ to $s_0+ 1$, we use (\ref{3.2}) to conclude that 
\begin{align}
\lim_{k\rightarrow \infty}\int_{s_0}^{s_0+ 1} \int_{M^n} s\Big|Rc_{g_{m_k}}+ Hess_{g_{m_k}} f_{m_k}- \frac{1}{2s} g_{m_k}\Big|^2 H_{m_k} d\mu_{g_{m_k}}(s)= 0 \nonumber
\end{align}

Therefore we have
\begin{align}
Rc_{g_{\infty}(s)}+ Hess_{g_{\infty}(s)} f_{\infty}= \frac{1}{2s} g_{\infty} \label{3.3}
\end{align}
for $s\in (s_0, s_0+ 1)$, where $f_{\infty}$ is defined by $H_{\infty}= (4\pi s)^{-\frac{n}{2}} e^{-f_{\infty}}$. Because $s_0$ is arbitrary positive number, (\ref{3.3}) holds for any $s\in (0, \infty)$. So $(M_{\infty}, g_{\infty}(s))_{s\in (0, \infty)}$ is a gradient shrinking soliton solution to backward Ricci flow.

Define $\ell_i(q, \theta)= \ell^{g_i}(q, \theta)$, where $\ell^{g_i}(q, \theta)$ is the reduced distance with base point $(y, 0)$ with respect to the backward Ricci flow solution $(M^n, g_i(t))_{t\in [0, \infty)}$. 

\begin{lemma}\label{lem l_1}
{If $(M^n, g(t))_{t\in [0, \infty)}$ is a non-flat Type I $\kappa$-solution to the backward Ricci flow for some $\kappa> 0$, and $Rm\geq 0$, then
\begin{align}
\ell_{\infty} (q, \theta)\vcentcolon = \lim_{i\rightarrow \infty}\ell_{i}(q, \theta) \label{l_1 1}
\end{align}
exists in the Cheeger-Gromov sense on $M_{\infty}\times(0, \infty)$, and $\ell_{\infty}(q, \theta)$ is the reduced distance of backward Ricci flow $(M_{\infty}, g_{\infty}(s))_{s\in (0, \infty)}$ with the base point $(y, 0)$.
}
\end{lemma}

\begin{remark}\label{rem 3.1}
{Note $(M_{\infty}, g_{\infty}(0))$ is the unique tangent cone at infinity of $(M^n, g(0))$, which is in fact an Euclidean metric cone (see Theorem $I.26$ in \cite{RFTA}). Although $(M_{\infty}, g_{\infty}(0))$ possibly is not a smooth manifold, from the definition of reduced distance and $(M_{\infty}, g_{\infty}(s))_{s\in (0, \infty)}$ are shrinking soliton solutions, the reduced distance of $(M_{\infty}, g_{\infty}(s))_{s\in (0, \infty)}$ with base point $(y, 0)$ is still well defined.
}
\end{remark}

\pf
{The similar argument as in the proof of Lemma $8.35$ and $8.36$ in \cite{RFG} leads to (\ref{l_1 1}), where $Rm\geq 0$ is required to get the uniform bounds of $\ell_i$ , $\nabla \ell_i$ and $\frac{\partial}{\partial \theta}\ell_i(q, \theta)$. Combining the definition of reduced distance and (\ref{l_1 1}), we get that $\ell_{\infty}(q, \theta)$ is the reduced distance of $(M_{\infty}, g_{\infty}(s))_{s\in (0, \infty)}$ with the base point $(y, 0)$.
}
\qed

In the rest of this section, $(M^n, g(t))_{t\in [0, \infty)}$ is a non-flat Type I $\kappa$-solution to the backward Ricci flow for some $\kappa> 0$, and $Rm\geq 0$.

\begin{lemma}\label{lem l_2}
{\begin{align}
|\nabla \ell_{\infty}(q, 1)|^2+ R_{g_{\infty}(1)}(q)= \ell_{\infty}(q, 1) \label{l_2 1}
\end{align}
and 
\begin{align}
\ell_{\infty}(q, 1)= f_{\infty}(q, 1)+ \beta \label{l_2 2}
\end{align}
where $\beta$ is a constant.
}
\end{lemma}

\pf
{By (\ref{3.3}) and Lemma \ref{lem l_1}, $\ell_{\infty}$ is the reduced distance on gradient shrinking soliton. From the argument in Lemma $7.77$ of Chapter $7$ in \cite{RFG} (especially $(7.113)$ and $(7.115)$ there), the conclusion follows (also see Section $3$ of \cite{CHI}).
}
\qed

\begin{prop}\label{prop 3.1}
{For $\beta$ in (\ref{l_2 2}), 
\begin{align}
\widehat{V}_{\infty}(y, 0)= e^{-\beta} \label{3.1.0}
\end{align}
}
\end{prop}

\begin{remark}\label{rem 3.2}
{From (\ref{3.1.0}), we know that $\beta$ is independent of the choice of sequence $\{t_i\}_{i= 1}^{\infty}$.
}
\end{remark}

\pf
{From (\ref{def V}) and Lemma \ref{lem l_1}, we can use Lemma $8.38$ of \cite{RFG} and (\ref{l_2 2}) to conclude 
\begin{align}
\widehat{V}_{\infty}(y, 0)&= \int_{M_{\infty}} (4\pi)^{-\frac{n}{2}} e^{-\ell_{\infty}(q, 1)} d\mu_{g_{\infty}(1)}(q) \nonumber \\
&= e^{-\beta} \int_{M_{\infty}} H_{\infty}(q, y, 1) d\mu_{g_{\infty}(1)}(q) \label{3.1.1}
\end{align}

From (\ref{2.4.9.1}) and the definitions of $H_i$, $g_i$, we obtain
\begin{align}
\int_{B_{g_i}(y, b, 1)} H_i(x, y, 1) d\mu_{g_i(1)}(x)\geq 1- C\exp{(-C b^2)} \label{3.1.2}
\end{align}
where $B_{g_i}(y, b, 1)= \{x| d_{g_i(1)} (x, y)< b,\ x\in M^n\}$, $b\geq 1$ and $C= C(C_1, n, \kappa)$.

Note we have
\begin{align}
\int_{B_{g_{\infty}}(y, 2b, 1)} H_{\infty} (x, y, 1) d\mu_{g_{\infty}(1)} (x) \geq \lim_{i\rightarrow \infty} \int_{B_{g_i}(y, b, 1)} H_i  \label{3.1.3}
\end{align}

From (\ref{3.1.2}) and (\ref{3.1.3}), 
\begin{align}
\int_{B_{g_{\infty}}(y, 2b, 1)} H_{\infty} \geq 1- C\exp{(-C b^2)} \nonumber
\end{align}
Let $b\rightarrow \infty$ in the above inequality, we get 
\begin{align}
\int_{M_{\infty}} H_{\infty} (x, y, 1) d\mu_{g_{\infty}(1)}(x) \geq 1 \nonumber
\end{align}

On the other side, we get $\int_{M_{\infty}} H_{\infty} (x, y, 1) d\mu_{g_{\infty}(1)}(x) \leq 1$ from Fatou's lemma. Hence 
\begin{align}
\int_{M_{\infty}} H_{\infty} (x, y, 1) d\mu_{g_{\infty}(1)}(x)= 1 \label{3.1.4}
\end{align}

From (\ref{3.1.1}) and (\ref{3.1.4}), we get our conclusion.
}
\qed

Define $\nu(t)= \Big[t(2\Delta f- |\nabla f|^2+ R)+ f- n\Big] H(t)$ and 
\begin{align}
\nu_{\infty}(x, y, 1)= \lim_{i\rightarrow \infty} \nu(t_i)t_i^{\frac{n}{2}}= \Big[2\Delta_{g_{\infty}(1)} f_{\infty}- |\nabla f_{\infty}|^2_{g_{\infty}(1)}+ R_{g_{\infty}(1)}+ f_{\infty}- n\Big] H_{\infty} (x, y, 1) \nonumber 
\end{align}
We do not know whether $\nu_{\infty}(x, y, 1)$ is independent on the choice of sequence $\{t_i\}$, but we have the following Lemma, whose original proof appeared in Section $3$ of \cite{CHI}.

\begin{lemma}\label{lem 3.2}
{\begin{align}
\int_{M_{\infty}} \nu_{\infty}(1) d\mu_{g_{\infty}(1)}= -\beta\ , \quad 
\int_{M_{\infty}} f_{\infty} H_{\infty}(1) d\mu_{g_{\infty}(1)} = \frac{n}{2}- \beta \label{3.2.2} 
\end{align}
}
\end{lemma}

\pf
{Do integration by parts (it is easy to justify integration by parts near infinity using the results in Section \ref{section 3}), then use $\Delta_{g_{\infty}(1)} f_{\infty}+ R_{\infty}= \frac{n}{2}$, we can get 
\begin{align}
\int_{M_{\infty}} \nu_{\infty}&= \int_{M_{\infty}} \Big[\Delta_{\infty} f_{\infty}+ R_{\infty}+ f_{\infty}- n\Big] H_{\infty} \nonumber \\
&= \int_{M_{\infty}} \Big[f_{\infty}- \frac{n}{2}\Big] H_{\infty} = \Big(\int_{M_{\infty}} f_{\infty} H_{\infty} \Big)- \frac{n}{2} \label{3.2.3}
\end{align}

On the other hand, do integration by parts, from (\ref{l_2 1}) and (\ref{l_2 2}), we get 
\begin{align}
\int_{M_{\infty}} \nu_{\infty} &= \int_{M_{\infty}} \Big[|\nabla f_{\infty}|^2+ R_{\infty}+ f_{\infty}- n\Big] H_{\infty} \nonumber \\
&= \int_{M_{\infty}} \Big[|\nabla \ell_{\infty}|^2+ R_{\infty}+ f_{\infty}- n\Big] H_{\infty} \nonumber \\
&= \int_{M_{\infty}} \Big[2f_{\infty}+ \beta- n\Big] H_{\infty}= 2\Big(\int_{M_{\infty}} f_{\infty}H_{\infty} \Big)+ (\beta- n) \label{3.2.4}
\end{align}
From (\ref{3.2.3}) and (\ref{3.2.4}), we get our conclusion.
}
\qed

From (\ref{W}), it is easy to check that 
\begin{align}
\mathscr{W}(g, f, t)= \Big(\int_{M^n} \nu \Big)(t)= t\frac{\partial}{\partial t}N(g, H, t)+ N(g, H, t) \label{3.4}
\end{align}

\begin{prop}\label{prop 3.3}
{\begin{align}
\lim_{t\rightarrow \infty} \mathscr{W}(g, f, t)= \lim_{t\rightarrow \infty} N(g, H, t)= -\beta \nonumber
\end{align}
}
\end{prop}

\pf
{We firstly show $\lim_{t\rightarrow \infty} N(g, H, t)= -\beta$.
\begin{align}
\Big|\int_{M^n\backslash B(y, b\sqrt{t}, t)} (fH)(x, y, t) d\mu_{g(t)}(x)\Big| &\leq \int_{M^n\backslash B(y, b\sqrt{t}, t)} \Big|\ln \big((4\pi t)^{\frac{n}{2}} H\big)\Big|\cdot H \nonumber \\
&\leq \int_{M^n\backslash B(y, b\sqrt{t}, t)} CH \leq C\exp{(-Cb^2)} \label{3.3.1}
\end{align}
where we used the bound of $t^{\frac{n}{2}}H$ in Lemma \ref{lem 4.1 of CZ} and (\ref{2.4.7}), $b\geq \frac{1}{2}$ is a constant and $C= C(C_1, n, \kappa)$ is independent of $b$ and $t$.

For any sequence $\{t_i\}$ as in the beginning of this section, we can find a subsequence also denoted as $\{t_i\}$, such that $g_{i}$ converges to $g_{\infty}$ in Cheeger-Gromov sense, and $H_i$ converges to $H_{\infty}$ on $M_{\infty}$ in $C^{\alpha}_{loc}$ topology.

Then using (\ref{3.3.1}) and the related convergence, we can get
\begin{align}
\lim_{i\rightarrow \infty} \Big(\int_{M^n} fH\Big)(t_i) &= \lim_{i\rightarrow \infty} \Big(\int_{M^n\backslash B(y, b\sqrt{t_i}, t_i)} fH\Big)(t_i)+ \lim_{i\rightarrow \infty} \Big(\int_{B(y, b\sqrt{t_i}, t_i)} fH\Big)(t_i) \nonumber \\
&\geq -C\exp\{(-Cb^2)\} + \lim_{i\rightarrow \infty} \Big[\int_{B_{g_i(1)}(b)} (fH)_{g_i(1)}\Big] \nonumber \\
&= -C\exp{(-Cb^2)}+ \int_{B_{g_{\infty}(1)}(b)} (f_{\infty} H_{\infty})_{g_{\infty}(1)} \nonumber
\end{align}
Let $b\rightarrow \infty$ in the above, we get
\begin{align}
\lim_{i\rightarrow \infty} \Big(\int_{M^n} fH\Big)(t_i)\geq \int_{M_{\infty}} f_{\infty}H_{\infty} (x, y, 1) d\mu_{g_{\infty}(1)}(x) \nonumber
\end{align}

On the other hand, similarly using (\ref{3.3.1}) and the convergence, we have
\begin{align}
\lim_{i\rightarrow \infty} \Big(\int_{M^n} fH\Big)(t_i)\leq \int_{M_{\infty}} f_{\infty}H_{\infty} (x, y, 1) d\mu_{g_{\infty}(1)}(x) \nonumber
\end{align}

By all the above and Lemma \ref{lem 3.2}, we get 
\begin{align}
\lim_{i\rightarrow \infty} \Big(\int_{M^n} fH\Big)(t_i)= \int_{M_{\infty}} f_{\infty}H_{\infty} (x, y, 1) d\mu_{g_{\infty}(1)}(x)= \frac{n}{2}- \beta \nonumber
\end{align}

From Proposition \ref{prop 3.1} we know that $\beta$ is independent of the choice of $\{t_i\}$, hence 
\begin{align}
\lim_{t\rightarrow \infty} \Big(\int_{M^n} fH\Big)(t)= \frac{n}{2}- \beta \nonumber
\end{align}
it is equivalent to
\begin{align}
\lim_{t\rightarrow \infty} N(g, H, t)=- \beta \label{3.3.2}
\end{align}

From (\ref{3.3.2}), we can get that $|N(g, H, 2t)- N(g, H, t)|\leq \epsilon$ for $t\gg 1$. This implies that there exists the sequence $\{t_i\}$ such that $t_i\frac{\partial}{\partial t}N(g, H, t_i)\rightarrow 0$ as $t_i\rightarrow \infty$. Hence from (\ref{3.4}),
\begin{align}
\lim_{t\rightarrow \infty} \mathscr{W}(g, f, t)= \lim_{i\rightarrow \infty} \mathscr{W}(g, f, t_i)= \lim_{i\rightarrow \infty}\Big[t_i\frac{\partial}{\partial t}N(g, H, t_i)+ N(g, H, t_i)\Big]= -\beta \nonumber
\end{align}
}
\qed

Combining Proposition \ref{prop 3.1} with Proposition \ref{prop 3.3}, Theorem \ref{thm 1.4} is proved.

\section*{Acknowledgments}
The author was partially supported by NSFC 11401336. He would like to thank Tobias H. Colding, Lei Ni for their comments, Bo Yang for his suggestion. He is also indebted to Ben Chow and Peter Li for discussion, Jiaping Wang for constant encouragement and support. He is particularly grateful to Xiaolong Li for his useful comments on the earlier version of the paper.

\end{document}